\documentclass[11pt]{article}
\usepackage{amsmath}
\usepackage{setspace}
\usepackage{fullpage}
\usepackage{graphicx}
\usepackage{amsthm}

\newtheorem{lemma}{Lemma}
\usepackage{epsf}
\newtheorem{thm}{Theorem}[section]
\newtheorem{prop}[thm]{Proposition}
\newtheorem{corr}[thm]{Corollary}
\newtheorem{q}[thm]{Question}
\theoremstyle{definition}

\begin{document}

\title{On low degree $k$-ordered graphs }

\author{
Karola M\'esz\'aros\\
Massachusetts  Institute of Technology\\
 {\tt karola@math.mit.edu}
\\
}

\date{}

\maketitle

%DELETEME start
%\LARGE
%DELETEME end

 \begin{abstract}
 
\begin{small}
A simple graph $G$ is \textit{k-ordered} (respectively, \textit{k-ordered
hamiltonian}) if, for any sequence of $k$ distinct vertices $v_1, \ldots, v_k$
of $G$, there exists a cycle (respectively, a hamiltonian cycle) in $G$
containing these $k$ vertices in the specified order.  In 1997 Ng and Schultz
  introduced these concepts of cycle orderability, and 
motivated by the fact that $k$-orderedness of a graph implies
$(k-1)$-connectivity, they
posed the
question of the existence of low degree $k$-ordered hamiltonian graphs.  
  We construct an infinite family of graphs,
 which we call
\textit{bracelet graphs}, that are $(k-1)$-regular and are $k$-ordered
hamiltonian for odd $k$. This result provides the best possible answer to the question of
the existence of low degree $k$-ordered hamiltonian graphs for odd $k$. We
further  show 
that for even $k$,
there exist no $k$-ordered bracelet graphs with minimum degree
 $k-1$ and maximum degree less than $k+2$, and we exhibit an infinite family of
bracelet graphs with minimum degree $k-1$ and maximum degree $k+2$ that are
$k$-ordered for even $k$.    A concept related to
$k$-orderedness, namely that of
\textit{k-edge-orderedness}, is likewise
strongly related to connectivity properties.
We study this relation  in both undirected and directed graphs,  and   give bounds
on the connectivity  necessary to imply $k$-(edge-)orderedness properties. 
  \end{small}
\end{abstract}

\section{Introduction}

The concept of $k$-ordered graphs  was introduced in 1997  by Ng and
Schultz \cite{ng}. A simple graph  $G$ is a graph without loops or multiple
edges, and it is called  \textit{hamiltonian} if there exists a cycle (called
a \textit{hamiltonian cycle}) that contains all vertices of $G$.   In this paper we
consider only finite simple graphs.  A simple graph $G$ is called
\textit{k-ordered} (respectively, \textit{k-ordered hamiltonian}) if, for any
sequence of $k$ distinct vertices $v_1, v_2, \ldots, v_k$ of $G$, there exists a
cycle (respectively, a hamiltonian cycle) in $G$ containing these $k$ vertices
in the specified order. Previous results concerning cycle orderability mainly
regard minimum degree and forbidden subgraph conditions that imply
$k$-orderedness or $k$-ordered hamiltonicity \cite{chen,  faun1, faun2,
faun-gould,  sar}.  A comprehensive survey of results can be found in \cite{faun}. 

A notion related to $k$-orderedness, that of  \textit{k-edge-orderedness},  
has been studied in \cite{c}.  A simple graph $G$ is
\textit{k-edge-ordered} (respectively, \textit{k-edge-ordered eulerian}) if,
for any sequence of $k$ distinct edges $e_1, e_2, \ldots, e_k$ of $G$, there
exists a tour (respectively, an eulerian tour, that is, a tour containing each
edge of $G$) in $G$ containing these $k$ edges in the specified order. 
It is natural to explore analogous notions in directed graphs. A directed
graph $D$ is \textit{k-ordered (hamiltonian)}   if, for any sequence of $k$
distinct vertices $v_1, v_2, \ldots, v_k$ of $D$, there exists a directed
(hamiltonian) cycle in $D$ containing these $k$ vertices in the specified
order. Furthermore,   $D$ is \textit{k-edge-ordered (eulerian)}   if, for any
sequence of $k$ distinct edges $e_1, e_2, \ldots, e_k$ of $D$, there exists a
directed (eulerian) tour in $D$ containing these $k$ edges in the specified
order.

As $k$-orderedness implies
$(k-1)$-connectivity, a natural question to pose is the existence of low
degree $k$-ordered graphs.  The question of the
existence of $3$-regular $4$-ordered graphs  was posed in \cite{ng} and
answered in the affirmative in \cite{karola}. In Section 2, we answer the more general question
of the existence of $(k-1)$-regular $k$-ordered graphs for odd $k$; in particular,
we exhibit an infinite family of graphs, called \textit{bracelet graphs}, that
are $(k-1)$-regular and $k$-ordered hamiltonian.  We also exhibit 
sufficient conditions for a bracelet graph to be $k$-ordered. 

In Section 3   we exhibit a bound on the diameter of a $k$-ordered graph, and we show 
that the bound is almost tight for the bracelet graph that we constructed in Section 2. 
In Section 4 we continue investigating low degree $k$-ordered graphs for even $k$, and 
we show that for
even $k$  there are no $k$-ordered bracelet graphs with minimum degree
$k-1$ and maximum degree less than $k+2$; however, we also exhibit an infinite family of
bracelet graphs with minimum degree $k-1$ and maximum degree $k+2$ that are
$k$-ordered for even $k$. This construction partially answers the question of the existence of
low degree $k$-ordered graphs for even $k$. 

In Section 5 we consider $k$-orderedness properties of directed graphs,
exhibiting an infinite family of $(k-1)$-diregular graphs that
are 
$k$-ordered hamiltonian. 
In Sections 6 and 7 
we establish a relation between connectivity and $k$-(edge-)orderedness in undirected as well as 
 directed graphs.  We conclude our paper by posing open questions.

\section{$2k$-regular $(2k+1)$-ordered hamiltonian graphs}

As observed in \cite{ng}, a $k$-ordered graph $G$ is also $(k-1)$-connected, and hence has 
minimum degree at least $k-1$.
The question of the existence of an infinite family of $3$-regular $4$-ordered
graphs was raised in \cite{ng} and answered in \cite{karola} by
constructing such a family. More generally, we are interested in whether there
exists an infinite family of $(k-1)$-regular $k$-ordered graphs. In this
section we answer this question in the case where $k$ is
odd,  exhibiting an infinite
family of $(k-1)$-regular $k$-ordered hamiltonian graphs for all odd $k\geq 3$.  

We call a graph $G$ a \textit{bracelet graph} if its  vertex set $V$ can be
partitioned into $V_1\cup V_2 \cup  \cdots \cup V_m$, $m\geq 3$, with  
$V_i$  nonempty for  all $i \in [m]$ (we denote the set $\{1, 2, \ldots, m\}$ by $[m]$),
  such that $v$
is adjacent to $u$  in $G$ if and only if   $v\in V_i$ and $u \in
V_j$  and $i-j \equiv 1$ or $-1$ (mod $m$).  We call $V_i$, for $i \in [m]$,
a \textit{part} of   $G$, and denote its cardinality by $|V_i|$.  We say that
two parts $V_i$ and $V_j$ are adjacent if $i-j \equiv 1$ or $-1$   (mod $m$). 
We also say that parts $V_i$ and $V_j$ are at distance $d$ if there
is a path from a vertex in $V_i$ to a vertex in $V_j$ such that it contains
$d$ edges and there are no $2$ vertices on the path from the same part. Note that
as bracelet graphs are ``cyclic'' there are two options for the distance between two
parts; in general, it will be clear from the context which of the two
distances we mean. 

Throughout this paper we will frequently want to construct a cycle or path through vertices in 
a specified order. We will refer to these specified vertices as \textit{marked vertices}. We also 
use the idea of free vertices in the course of the paper; we shall define free vertices 
in the statement 
of Lemma~\ref{free}, which we will use for proving Theorem~\ref{2.1}. 

Let $G_{k, 2m}$ be a bracelet graph with parts $V_1, V_2, \ldots,
V_{2m}$, $m\geq 2$, such that $|V_i|=k$ for $i \in [2m]$.  It is
clear that $G$ is simple and  $2k$-regular by construction.

\begin{lemma} \label{free}
  Given $2k+1$ marked vertices $v_1$,
$v_2$, \ldots, $v_{2k+1}$ in
$G_{k, 2m}$,  there exists a set of $2m$ 
vertices, which we call \emph{free vertices}, satisfying the following two
properties:  (i) there is exactly one free vertex in
each part of $G_{k, 2m}$, and (ii) there exists some $i$ such that the marked
vertices $v_i$ and $v_{i+1}$ are in the set of free vertices (indices taken
modulo $2k+1$) and no other marked vertices are in the set of free vertices.
\end{lemma}

\begin{proof}
 If each part of $G_{k, 2m}$ contains at most $k-1$
marked vertices let $B_1$  any part that contains a marked vertex, and 
 if there is exactly one part that has $k$ marked vertices, let $B_1$ be that part. 
Take
a vertex $v_i$ from $B_1$ such that $v_{i+1} \notin B_1$. Then we can take
vertices $v_i, v_{i+1}$, and one unmarked vertex from each of the parts
not containing $v_i$ or $v_{i+1}$ as the free vertices. 

If there are exactly two parts, $B_1$ and $B_2$,  that have $k$ marked
vertices, then it is not hard to see that there are two consecutive vertices
$v_i$ and $v_{i+1}$  (indices taken modulo $2k+1$) such that $v_i \in B_1$ and
$v_{i+1} \in B_2$.  Then we can take vertices $v_i, v_{i+1}$, and one
unmarked vertex from each of the parts not containing $v_i$ or $v_{i+1}$  as
the free vertices. This completes the proof.
\end{proof}

\begin{thm} \label{2.1}
For every $k\geq 1$, there exists an infinite family of $2k$-regular graphs
that are $(2k+1)$-ordered hamiltonian.  \end{thm}

\begin{proof} We prove that the bracelet graphs $G_{k, 2m}$ introduced above 
are  $(2k+1)$-ordered hamiltonian  for $k
\geq 1$ and $m \geq 2$.  In fact, we will prove more: given any $2k+1$ vertices
$v_1, v_2, \ldots, v_{2k+1}$ in $G_{k, 2m}$, there exists a hamiltonian cycle
$\mathcal{H}_{k, 2m}$ of $G_{k, 2m}$ that traverses the vertices in order, and 
satisfies the following condition ($\star$): 
 for any two adjacent parts $B_1$ and $B_2$ of $G_{k, 2m}$, there exists
an edge of $\mathcal{H}_{k, 2m}$ with one vertex in each of $B_1$ and $B_2$. 
We will proceed by  induction on $k$. 

\textit{Base case: $k=1$.} 

It is clear that $G_{1, 2m}$ for $m\geq 2$ is just a cycle, and it follows
that $G_{1,2m}$ is $3$-ordered hamiltonian; futhermore we can
take $\mathcal{H}_{1, 2m}=G_{1, 2m}$. 

\textit{Inductive step.} 

Suppose that $G_{k-1, 2m}$ is $(2k-1)$-ordered hamiltonian for $m\geq 2$ and
given $2k-1$ vertices in $G_{k-1, 2m}$ there is a  hamiltonian cycle
$\mathcal{H}_{k-1, 2m}$ satisfying condition ($\star$). 
Consider   the $2k+1$ marked vertices $v_1, v_2, \ldots, v_{2k+1}$ through which we 
wish to construct a hamiltonian cycle. By Lemma~\ref{free}
it is possible to find $2m$ free vertices in $G_{k, 2m}$, one in each part. Without loss of
generality we can suppose that the two marked vertices among the $2m$ free
vertices are $v_{2k}$  and $v_{2k+1}$. Note that the  graph induced by the
$2m$ free vertices is a cycle, $\mathcal{C}$, and the graph   induced on the
vertices of $G_{k, 2m}$ without the free vertices is (isomorphic to) $G_{k-1,
2m}$. Therefore, by the induction hypothesis,
 there exists a hamiltonian cycle $\mathcal{H}_{k-1, 2m}$
through the $2k-1$ marked vertices $v_1, v_2, \ldots, v_{2k-1}$ in $G_{k-1,
2m}$ satisfying condition ($\star$).
 
 We show that
given  $v_1, v_2, \ldots, v_{2k+1}$ in $G_{k, 2m}$ there is a hamiltonian cycle
$\mathcal{H}_{k, 2m}$ containing the $2k+1$ vertices in the specified order,
and such that for any two adjacent parts $B_1$ and $B_2$ of $G_{k, 2m}$ there
exists an edge of $\mathcal{H}_{k, 2m}$ with one vertex in each of $B_1$ and
$B_2$.  This will also show, in particular, that $G_{k, 2m}$ is
$(2k+1)$-ordered hamiltonian for $m\geq 2$.  We will examine 
cases depending on the positions of the $2k+1$ specified vertices and show how to
construct the desired hamiltonian cycle in each case.
 
\textit{Case 1.} Suppose $v_{2k-1}$ and  $v_{2k}$ are in different parts. In
this case the hamiltonian cycle  $\mathcal{H}_{k, 2m}$ in $G_{k, 2m}$ is as
follows.  Follow $\mathcal{H}_{k-1, 2m}$ in $G_{k-1, 2m}$ from $v_1$
until reaching
$v_{2k-1}$. 
If $v_{2k}$ is in a part adjacent to the part of $v_{2k-1}$, go to $v_{2k}$ from
$v_{2k-1}$ and   continue going to
the free vertices in the not yet  visited adjacent parts along the cycle
$\mathcal{C}$ so that we reach $v_{2k+1}$. 
If $v_{2k}$ is not in a part adjacent to the part containing $v_{2k-1}$, then  continue going to
the free vertices in the not yet  visited adjacent parts along the cycle
$\mathcal{C}$ so that we reach $v_{2k}$ first and then $v_{2k+1}$. 
In both cases (whether or not $v_{2k}$ is in a part adjacent to the part containing $v_{2k-1}$)
continue along $\mathcal{C}$  after meeting $v_{2k+1}$ until reaching the
free vertex of   the part containing $v_{2k-1}$.  After this, go to the
vertex that is adjacent to $v_{2k-1}$  in $\mathcal{H}_{k-1, 2m}$ when
going from $v_{2k-1}$ to $v_1$, and continue on $\mathcal{H}_{k-1, 2m}$ until
$v_1$.  It is clear that the hamiltonian cycle $\mathcal{H}_{k, 2m}$ that we
have constructed has the property that for any two adjacent parts of $G_{k, 2m}$ there exists
an edge of $\mathcal{H}_{k, 2m}$ with one vertex in each of the two adjacent
parts.

\textit{Case 2.} Suppose $v_{2k-1}$ and $v_{2k}$ are in the same part.

Let the part containing $v_{2k-1}$ and $v_{2k}$ be $B_1$, and the part
containing $v_{2k+1}$ be $B_2$. Note that $B_1\neq B_2$ since $v_{2k}$ and 
$v_{2k+1}$ are in different parts. Let $u$ be the vertex adjacent to 
  $v_{2k-1}$ in the traversal of   $\mathcal{H}_{k-1, 2m}$ from
$v_{2k-1}$ to $v_1$.

\textit{Case 2.1.} Suppose $u \not \in B_2$.  Follow $\mathcal{H}_{k-1, 2m}$
 in $G_{k-1, 2m}$ from $v_1$ to $v_{2k-1}$. From $v_{2k-1}$ go to an unmarked
 free vertex in the part containing $u$.  From this free vertex go to $v_{2k}$, and after this
 continue going to the free vertices in the not yet  visited adjacent parts
 along $\mathcal{C}$. During this we meet $v_{2k+1}$, and continue along
 $\mathcal{C}$ until we are at the free vertex of   the part adjacent to the
 part containing $u$.  Go to $u$ and continue on $\mathcal{H}_{k-1, 2m}$ until $v_1$.
 The obtained hamiltonian cycle is $\mathcal{H}_{k, 2m}$.

\textit{Case 2.2.} Suppose $u \in B_2$.  
In this case, $B_1$ and $B_2$ are adjacent.  Follow $\mathcal{H}_{k-1, 2m}$ in
$G_{k-1, 2m}$ from $v_1$ to $v_{2k-1}$. From $v_{2k-1}$ go to an unmarked free
vertex in the adjacent part that is not $B_2$. From  here go to $v_{2k}$,
then to $v_{2k+1}$, then to the free vertex in the next adjacent part and back to
$u$. Continue on $\mathcal{H}_{k-1, 2m}$ until reaching $v_1$.  When $m=2$ this is a
hamiltonian cycle in $G_{k, 2m}$, but if $m>2$, this is not a hamiltonian
cycle.  We now show that  we can reroute the path from $v_1$ to $v_{2k-1}$ so
that we pick up all the missing vertices.

Indeed, note the following. If there is an edge $ab$ in the hamiltonian cycle
$\mathcal{H}_{k-1, 2m}$ in $G_{k-1, 2m}$ and $c$ and $d$ are free unmarked
vertices, such that $a$ and $c$ are in the same part and $b$ and $d$ are in
the same part, then replacing edge $ab$ by edges $ad$, $dc$, and $cb$, preserves
the ordering, and includes $c$ and $d$ in the cycle. Call this operation of
rerouting  $\alpha$. 

By the inductive hypothesis,  $\mathcal{H}_{k-1, 2m}$ in $G_{k-1, 2m}$  is such
that for any two adjacent parts $B_1$ and $B_2$ in $G_{k-1, 2m}$ there exists an
edge of $\mathcal{H}_{k-1, 2m}$ with one vertex in each of $B_1$ and $B_2$.
Because the number of parts is even, we can pair up
adjacent parts
(without the part containing  $v_{2k}$, the part containing $v_{2k+1}$ and the two parts adjacent
to these) and perform the rerouting 
operation $\alpha$ as explained in the preceding paragraph for the
$m-2$ part pairs.  The hamiltonian path $\mathcal{H}_{k, 2m}$ satisfies
condition ($\star$), concluding the proof.

\end{proof}

\begin{corr}
The graphs $G_{k, m}$ are $(2k+1)$-ordered for all $m \geq 4$.
\end{corr}

\begin{proof}
Note that in the proof of  Theorem~\ref{2.1} the assumption about the even number of 
parts was used only in  Case 2.2 where we  rerouted the cycle using $\alpha$.
To prove $(2k+1)$-orderedness, we do not require the rerouting, thus 
$G_{k, m}$ is $(2k+1)$-ordered for all $m \geq 4$.
 \end{proof}

\begin{thm}\label{2.3}
Any bracelet graph with at least $4$ parts, at least $k$ vertices in
each part, and at least $2k+1$ total vertices in every pair of parts at distance 2 is
$(2k+1)$-ordered.
\end{thm}

\begin{proof} Let $G$ be a bracelet graph with parts $V_1, V_2, \ldots, V_m$, such
 that each part has at least $k$ vertices and there are at least $2k+1$
 vertices in any two parts at distance $2$.  Let $v_1, v_2, \ldots, v_{2k+1}$ be
 any $2k+1$ specified vertices. We show that there exists a cycle in $G$
 containing $v_1, v_2, \ldots, v_{2k+1}$ in this order, and therefore  
 $G$ is $(2k+1)$-ordered. 

We will prove the statement by induction on $k$. 
 
\textit{Base case: $k=1$.} It is not hard to see that any bracelet graph
 with at least $1$ vertex per part and at least $3$ vertices in any two
 parts at distance $2$ is $3$-ordered. 

\textit{Inductive step.} 

\textit{Case 1.} Suppose each part contains at most $k$ marked vertices.  Let
$V'_i \subset V_i$ for all $i \in [m]$ such that each $V'_i$ contains $k$ vertices,
and all of the marked vertices are in $V'=V'_1 \cup V'_2 \cup \cdots \cup
V'_m$.  As the
graph induced by $V'$ is isomorphic to $G_{k, m}$, it follows by Corollary 2.3
that we can find the desired cycle. 

\textit{Case 2.} Suppose there is a part $B$ that contains $k+l$ marked
vertices, $l\geq 1$.  As there are a total of $2k+1$ marked vertices, there can be
only one such part. 

\textit{Case 2.1.} Suppose  $B$ contains $2k+1$ marked vertices. As the two
 adjacent parts contain at least $2k+1$ vertices, it is not hard to see
 that there is a cycle containing the $2k+1$ vertices in the specified
 order. 

\textit{Case 2.2.} Suppose  $B$ contains fewer than $2k+1$ marked vertices.
 Then there exists a different part $B'$ such that $v_i \in B$ and $v_{i+1} \in B'$
 (indices taken modulo $2k+1$) for some $i$. Without loss of generality $v_{2k} \in B$
 and $v_{2k+1} \in B'$. In all of the other parts choose one vertex that is
 not marked (note that this is possible as all of the 
other parts contain at most $k-1$ marked vertices), and call it a free vertex. Consider the  graph $H$ induced on the
 vertices of $G$  without $v_{2k}, v_{2k+1}$, and the free vertices. This graph
 is $(2k-1)$-ordered by the inductive hypothesis.  Let $\mathcal{D}$ be the
 cycle containing vertices $v_1, v_2, \ldots, v_{2k-1}$ in this order in $H$.  Let
 $\mathcal{C}$ be the cycle induced on $v_{2k},
 v_{2k+1}$, and the free vertices. 

Now we show how to construct a cycle in $G$ containing the vertices  $v_1, v_2, \ldots,
v_{2k+1}$ in this order.

\textit{Case 2.2.1.} Suppose $v_{2k-1}$ and  $v_{2k}$ are in
 different parts.  In this case the cycle in $G$ is as follows. Follow
 $\mathcal{D}$ in $H$ from $v_1$ to $v_{2k-1}$. From $v_{2k-1}$ go to a free
 vertex $f$ in one of the adjacent parts, choosing the adjacent part so that
 either $f = v_{2k}$ or $f$ is
 unmarked. If $f\neq v_{2k}$, then  continue going to the free vertices in the not yet
 visited adjacent parts along the cycle $\mathcal{C}$ so that we
 reach
 $v_{2k}$ first and then $v_{2k+1}$. Continue along $\mathcal{C}$  after
 meeting $v_{2k+1}$ until we reach the free vertex of  the part containing
 $v_{2k-1}$.  After this, go to the vertex  that is adjacent to $v_{2k-1}$ 
in $\mathcal D$ when going
 from $v_{2k-1}$ to $v_1$, and continue along
 $\mathcal{D}$ until $v_1$.

\textit{Case 2.2.2.} Suppose $v_{2k-1}$ and $v_{2k}$ are in the same part. 
Denote the part containing $v_{2k-1}$ and $v_{2k}$ by $B_1$, and denote the part
containing $v_{2k+1}$ by $B_2$. Let $u$ be the end vertex of the edge incident
to $v_{2k-1}$ that is an edge in the segment of $\mathcal{D}$ from $v_{2k-1}$
to $v_1$.

\textit{Case 2.2.2.1} Suppose $u \not \in B_2$.  Follow $\mathcal{D}$ in $H$
from $v_1$ to $v_{2k-1}$. From $v_{2k-1}$ go to an unmarked free vertex in the
part containing $u$.   From the unmarked free vertex in
the part containing $u$ go to $v_{2k}$, and after this  continue going to
the free vertices in the not yet  visited adjacent parts along $\mathcal{C}$.
During this we meet $v_{2k+1}$, and we continue along $\mathcal{C}$  until we reach
the free vertex of   the part adjacent to the part of $u$.  Go to $u$ and
continue on $\mathcal{D}$ until $v_1$. 

\textit{Case 2.2.2.2.} Suppose $u \in B_2$.  
Then $B_1$ and $B_2$ are adjacent.  Follow $\mathcal{D}$ in $H$ from
$v_1$ to $v_{2k-1}$. From $v_{2k-1}$ go to an unmarked free vertex in the
adjacent part that is not $B_2$. From  here go to $v_{2k}$, then to
$v_{2k+1}$, then to the free vertex in the adjacent part and back to $u$.
Continue on $\mathcal{D}$ until $v_1$.  \end{proof}

\begin{corr}\label{a-cor-name}
 Any bracelet graph with at least $4$ parts and at least $k+1$ vertices in each part 
 is $(2k+1)$-ordered.
\end{corr}

In Corollary 2.5 we cannot  replace the requirement of each part
having at least $k+1$ vertices with each part having at least $k$ vertices. 
For a a subset $S$ of vertices let $N(S)$ denote 
the set of vertices adjacent to some vertex in  $S$, not including 
  vertices in $S$.  
If $v_1, v_2, \ldots, v_{2k+1}$ are independent vertices (there are no edges 
between any of them), and  if there 
is a cycle containing $v_1, v_2, \ldots, v_{2k+1}$ in this order 
 then $|N(\{v_1, v_2, \ldots, v_{2k+1}\})|\geq 2k+1$.  
Consider the example of a bracelet graph $G$ with  three adjacent parts $B_1$, $B_2$, and $B_3$
such that $|B_1|=|B_3|=k$ and $|B_2|=2k+1$. Specify  $v_1, v_2, \ldots,
v_{2k+1}$ to be in $B_2$. Then $v_1, v_2, \ldots, v_{2k+1}$ are independent, and 
then $|N(\{v_1, v_2, \ldots, v_{2k+1}\})|=2k$, which shows that such a bracelet graph cannot be 
$(2k+1)$-ordered. 
We also cannot weaken the condition in Theorem~\ref{2.3}  that
every two parts at distance $2$ have at least $2k+1$ vertices total. 

\begin{lemma} \label{ev} Given two vertices $u$ and $v$ in a bracelet graph $G$
with an even number of parts, the parity of the number of vertices on any
path from $u$ to $v$ is the same. Moreover,  any cycle in a 
bracelet graph $G$
with an even number of parts has even length. 
\end{lemma}

Although Lemma~\ref{ev} is an easy observation, it will  be a
convenient tool for proving non-existence of 
bracelet graphs with certain properties.

\begin{corr} There is no $k$-ordered hamiltonian bracelet graph $G$ with 
an even number of parts and an odd number of vertices.  \end{corr}

\section{The diameter of a $k$-ordered graph}

 In this section we give an upper bound on the diameter of
	$k$-ordered graphs, and we study the tightness of this bound. It has been
	observed by Denis Chebikin (personal communication)  that if $n$ is the number
	of vertices in a $4$-ordered graph, then the diameter of the graph is at most
	$\frac{n}{4}+2$. Indeed, as shown in the following proposition and corollary, a
	$k$-ordered graph has diameter not more than roughly $\frac{n}{k}$.

\begin{prop}\label{diam}
Given a $2k$-ordered graph $G$ on $n$ vertices, the diameter $d$ of $G$ is at most 
 $\lfloor \frac{n-3}{2k} \rfloor+2$.
\end{prop}

\begin{proof} Take two vertices $a$ and $b$ of $G$ at distance $d$. As $G$ is
$2k$-ordered, it is $(2k-1)$-connected, so in particular  both $a$ and $b$ have at
least $2k-1$ neighbors. Let $a_1, a_2, \ldots, a_{k-1}$ be $k-1$ distinct neighbors
of $a$, and let $b_1, b_2,  \ldots, b_{k-1}$ be $k-1$ distinct neighbors of $b$.
Note that the distance between $a_i$ and $b$ is at least $d-1$ for all $i\in
[k-1]$, and likewise the distance between $b_i$ and $a$ is at least $d-1$ for all
$i\in [k-1]$. Furthermore, the distance between  $a_i$ and $b_j$ is at least $d-2$
for all $i, j \in [k-1]$. Consider any cycle 
containing  the $2k$ vertices 
	$a, b_1, a_1, b_2, a_2, \ldots, b_{k-1}, a_{k-1}, b$ in the specified order. 
It contains at least $2k+(d-2)+(2k-3)(d-3)+(d-2)+(d-1)$
	vertices, so $n\geq 2kd-4k+4$, and $d\leq \frac{n}{2k}+2-\frac{4}{2k}$.
	The result follows. 

\end{proof}

\begin{corr}
Given a $(2k+1)$-ordered graph $G$ on $n$ vertices, the diameter $d$ of $G$ does 
not exceed $\lfloor \frac{n-3}{2k}\rfloor+2$.
\end{corr}

\begin{proof}
The statement follows from Proposition~\ref{diam}, as a $(2k+1)$-ordered graph is $2k$-ordered as well.
\end{proof}

\begin{thm} 
There exist infinitely many $(2k+1)$-ordered graphs on $n$ vertices 
with diameter  $\lfloor \frac{n-3}{2k}\rfloor+1$ for all $k\geq 2$.
\end{thm}

\begin{proof} Consider graphs $G_{k, 2m}$ as constructed in Section
2. In this
case $n=2mk$, and it is not hard to see that the diameter of $G_{k, 2m}$ is
$m=\lfloor\frac{2mk-3}{2k}\rfloor+1$.
\end{proof}

\section{Low degree $2k$-ordered graphs}

In this section we focus on low degree $2k$-ordered graphs.
Since $2k$-orderedness of a graph $G$ implies $(2k-1)$-connectivity, and thus
a $2k$-ordered graph $G$ has
minimum degree at least $2k-1$, a question of interest is the existence of 
$(2k-1)$-regular $2k$-ordered graphs, or low degree $2k$-ordered
graphs in general. Since $(2k+1)$-ordered  graphs are $2k$-ordered as well, it follows from
Theorem~\ref{2.1} that there is an infinite family of $2k$-regular graphs that are
$2k$-ordered. In this section we present a stronger result (Theorem~\ref{3.7}),
exhibiting an infinite  family of $2k$-ordered bracelet graphs such that the
minimum degree of these graphs is $2k-1$, and the maximum degree is $2k+2$. 
We also show that this result is the best possible for
bracelet graphs.

Note that the construction in Theorem~\ref{2.1} is specific to $2k$-regular graphs.
For example, consider an analogue of the graphs $G_{k, m}$ that are $(2k-1)$-regular. 
Note that in this case the number of vertices in all the parts cannot be the same. 
Thus, consider the graph $H_{k, 4m}$ to be a graph with $4m$ parts having a repeating pattern
of two parts with $k-1$ vertices followed by two parts with $k$ vertices.
Note that $H_{k, 4}$ is the bipartite graph $K_{2k-1, 2k-1}$, which was shown
in \cite{ng} to be $2k$-ordered hamiltonian. However, if $m>1$, then $H_{2,
4m}$ is not $4$-ordered, as it contains a square ($4$-cycle), and by \cite{karola} a
$3$-regular $4$-ordered graph on more than $6$ vertices cannot contain a
square.  In fact, it is easy to see that $H_{k, 4m}$ is not $2k$-ordered for
any $m>1$, because $2k$-orderedness implies $(2k-1)$-connectivity and the
deletion of two non-adjacent parts of $k-1$ vertices would disconnect $H_{k,
4m}$.

\begin{lemma} \label{eve} Let $G$ be a $2k$-ordered graph. If $s \leq k$, then there
exists no subset $V_1$ of the vertices of $G$ such that  $|V_1|=s$, 
$|V \setminus (N(V_1) \cup V_1)| \geq s$, and
$|N(V_1)|<2s$.  \end{lemma}

\begin{proof} Suppose  that $G$ is a graph such that there is a subset $V_1$
of the vertices of $G$ with $|V_1|=s$, 
$|V \setminus (N(V_1) \cup V_1)| \geq s$, and
$|N(V_1)|<2s$.    Specify $2s$ vertices $v_1, \ldots, v_{2s}$
such that for odd $i$, the $v_i$ are in
$V_1$, and for even $i$, the $v_i$ are
in $V \setminus (N(V_1)\cup V_1)$. Then it would be impossible to have a cycle containing
them in the specified order, because each vertex in $V_1$ would have to be adjacent to
$2$ distinct vertices in $N(V_1)$. Thus, $G$ could not have been $2s$-ordered,
nor $2k$-ordered for $k \geq s$.  \end{proof}

\begin{corr}  \label{3.2} Given a $2k$-ordered bracelet graph $G$ with more than $5$ parts,
there exists no part $B$ in $G$ such that $|B| \leq k$ and $N(B) < 2|B|$, and
there also exists no part $B'$ in $G$ such that $|B'| > k$ and $N(B') < 2k$.
\end{corr}

 \begin{proof} 
 Consider $6$ consecutive parts $B_1, B_2, \ldots, B_6$ of a $2k$-ordered 
bracelet graph $G$. Since $2k$-orderedness implies $(2k-1)$-connectivity, 
 the total number of vertices in parts $B_4$ and $B_6$ is at least $2k-1$.

Suppose that 
 there exists a part $B$ in $G$ such that $|B| \leq k$ and $N(B) < 2|B|$, and let $B=B_2$ 
 without loss of generality.  Taking 
$V_1$ as described in Lemma~\ref{eve} to be the vertices in $B_2$ it 
 follows that $G$ cannot be $2k$-ordered, contradicting our assumption. 

Suppose that there exists a part $B'$ in $G$ such that $|B'| > k$ and $N(B') < 2k$, and 
let $B'=B_2$, without loss of generality. 
Taking 
$V_1$ as described in Lemma~\ref{eve} to be some $k$  vertices in $B_2$ it 
 follows that $G$ cannot be $2k$-ordered, contradicting our assumption. 
 
 \end{proof}

\begin{prop} \label{3.3} There is no $2k$-ordered bracelet graph with more than 6 parts that has minimum
degree $2k-1$ and maximum  degree less than $2k+2$. 
\end{prop}

\begin{proof} Suppose that there exists a $2k$-ordered bracelet graph $G$ with more than 6
parts that has minimum degree $2k-1$ and maximum  degree less than $2k+2$.

\textit{Claim.} There is no part containing fewer than $k-1$ vertices. 

Let $B$ be a part with the minimum number of vertices, $b$, in bracelet graph $G$. 
Let $B$, $B_1$, $B_2$, $B_3$, and $B_4$ be a sequence of $5$
adjacent parts. Then
$B_2$ contains $2k+1-b$, $2k-b$ or $2k-1-b$ vertices, while $B_4$ contains
$b$, $b+1$ or $b+2$ vertices (since $B$ was  a part with the minimum number
of vertices).  Thus, because $G$ is $(2k-1)$-connected and there are at least $6$ parts, 
 $b+(b+2)\geq 2k-1$,
implying that $b \geq k-1$.    

Therefore, since there exists a vertex in $G$ with degree $2k-1$, there exist parts $B'$ and $B'_2$ 
distance $2$ apart, with $k-1$ and $k$ vertices respectively.   
Applying Corollary~\ref{3.2} to $B'_1$, the part between $B'$ and $B'_2$ , we see that 
$B'_1$ also has $k-1$ vertices. 
Also, since the graph is $(2k-1)$-connected, every part  other than 
$B'$ and $B'_1$ must have at least $k$ vertices.

Suppose part $C$ contains
$k$ vertices, and $C$ is not $B'_2$ or the part adjacent to $B'$ (that is not $B'_1$). 
Let $C_1$ be the part following $C$, in the same direction that 
$B'_1$ follows $B'$.  Choose 
$v_2, v_4, \ldots, v_{2k-2}$ in $B'$, choose $v_1, v_3, \ldots,
v_{2k-1}$  in $B'_2$, and choose $v_{2k}$ in $C_1$. Then it is not hard to see that
there can be no cycle visiting the $2k$ vertices in order, 
as there could be no path leading from $v_{2k}$ to
$v_1$ once we have traversed $v_1, v_2, \ldots, v_{2k}$. 
Thus, the only parts that might have fewer than $k+1$ vertices are $B', B'_1, B'_2$ and the parts 
adjacent to $B'$. Since there are at least $7$ parts,   it follows that
there is a vertex with degree at least $2k+2$, contradicting our assumptions.  
\end{proof}

By arguments analogous to those in Proposition~\ref{3.3}, it is not hard to see
 that there is no $(2k-1)$-regular $2k$-ordered bracelet graph on more than $4$ parts, 
 and that there is no $2k$-ordered bracelet graph with minimum degree $2k-1$ and maximum degree
 $2k$ on more than $5$ parts. Also, the $k$-orderedness of 
bracelet graphs on at most $6$ parts can be easily studied.  In the following theorem we show that 
 Proposition~\ref{3.3} is the best possible by exhibiting an infinite family of 
 $2k$-ordered graphs with minimum degree $2k-1$ and maximum degree
 $2k+2$.

\begin{thm} \label{3.7}
There exists an infinite family of $2k$-ordered graphs $P_{k, m}$ with
   minimum degree $2k-1$ and maximum
 degree $2k+2$.
 \end{thm}

\begin{proof} 
For $m\geq 5$, let $P_{k, m}$ be the bracelet graph with $m$ parts, $V_1, V_2, \ldots, V_m$, 
satisfying 
 $|V_1|=|V_2|=k-1$, $|V_3|=k$, and $|V_i|=k+1$ for all $i>3$. 
 We will show by induction on $k$ that 
the graphs $P_{k, m}$ are $2k$-ordered for $k \geq 2$ and $m \geq 5$.

\textit{Base case: $k=2$.} 

\textit{Case 1.} Suppose one of the parts contains $3$ of the $4$ marked vertices. Call this
part $B$.  Without loss of generality we may assume that it contains vertices  $v_1, v_2$, and
 $v_3$.
As there is a part adjacent to $B$ with $3$ vertices, it is possible to go
from $v_1$ to $v_3$ through $v_2$,  meeting two vertices from that adjacent
part, and then it is possible to go from $v_3$ by going through all  the parts
to $v_4$ and on to $v_1$. 

\textit{Case 2.} Suppose one of the parts contains $2$ marked vertices. Call such a
part $B$.  Without loss of generality we may assume that it contains vertices $v_1$ and $v_2$ or
$v_1$ and $v_3$. 

\textit{Case 2.1.} Suppose $B$ contains vertices $v_1$ and  $v_2$. 

Using one vertex (that is not $v_3$ or $v_4$)
from an adjacent part with $3$ vertices we can go from $v_1$ to $v_2$, and
then,  regardless of where $v_3$ and $v_4$ are, it is not hard to see that the
desired cycle exists.

\textit{Case 2.2.} Suppose $B$ contains vertices $v_1$ and  $v_3$. We can divide this case up 
into cases depending on whether 
  $v_2$ or $v_4$ are in parts with $1$ or $3$ 
vertices, or whether they are both in the part with $2$ vertices. We can easily find the desired 
cycles in each case.

\textit{Case 3.} Suppose all four vertices are in different parts.  Let
$v_1 \in V_{i_1}$, $v_2 \in V_{i_2}$, $v_3 \in V_{i_3}$, and $v_4 \in V_{i_4}$. 
Without loss of generality we can suppose that $0<i_2-i_1<i_3-i_1<i_4-i_1$ 
or $0<i_2-i_1<i_4-i_1<i_3-i_1$ (as by symmetry we can rotate and reflect an ordered cycle), 
where we consider subtraction modulo $m$ taking results between $0$ and $m-1$. 
In the case $0<i_2-i_1<i_3-i_1<i_4-i_1$ it is clear that there is a cycle containing the four
vertices in the specified order, and if  $0<i_2-i_1<i_4-i_1<i_3-i_1$, 
analysis  shows the existence of the desired cycle.

\textit{Inductive step.} Suppose the claim is true for all numbers less than
$k$. We shall show that $P_{k, m}$ is still $2k$-ordered for $m \geq 5$.  

Suppose there is no part with all vertices marked or there is exactly one
part with all vertices marked. Then we can find two different parts
containing two consecutive vertices, $v_{2k-1}$ and $v_{2k}$ without loss of
generality, such that all other parts have a free vertex.  By arguments
analogous to those in Theorems~\ref{2.1} and ~\ref{2.3}  we can show how to find a cycle in
$P_{k, m}$ using the cycle in $P_{k-1, m}$.  

If there are exactly two parts  with all vertices marked, then 
one of the following cases occurs:

 (i)  both parts have exactly $k-1$ vertices,

(ii)  one  part has exactly $k-1$ vertices and the other has exactly $k$ vertices, 
or 

(iii) 
 one of the  parts has exactly $k-1$ vertices and the other has exactly $k+1$ vertices.

If either case (ii) or
case (iii) occurs, then there exist marked vertices $v_i$ and $v_{i+1}$, one in 
either part with all vertices marked.
The result then follows from arguments similar to those in
Theorems~\ref{2.1} and \ref{2.3}. If,
however, the two parts  with all vertices marked both have $k-1$
vertices
and there are no two vertices  in them adjacent (if there are, then we are
done by arguments as before), then, without loss of generality,  the two parts
with all vertices marked contain the vertices $v_2, v_3, \ldots, v_k$ and $v_{k+2}, v_{k+3}, \ldots, v_{2k}$
respectively. It is not hard to see that in this case there exists a cycle
containing $v_1, v_2, \ldots, v_{2k}$ in this order, regardless of where $v_1$
and $v_{k+1}$ are. 

\end{proof}

\section{$k$-Ordered directed graphs}
 
In this section we address    the
existence of low degree $k$-ordered directed graphs.
 This inquiry is motivated by the analogous
questions for undirected graphs. 

Consider a directed graph $D$, and denote its set of vertices by $V(D)$.
A directed graph $D$ is said to be \textit{strongly connected} 
if given any two vertices 
 $u$ and $v$ in $D$, there exists a directed path from $u$ to $v$.
A \textit{vertex cut} of a digraph $D$ is a set $S\subset V(D)$  
  such that 
$D-S$ is not strongly  connected. 
 
\begin{prop}
If a directed graph $D$ is $k$-ordered, then every vertex cut has at least $k-1$ vertices.  
\end{prop}

The proof is analogous to the proof of the undirected case in \cite{ng}.

\begin{corr} If a directed graph $D$ is $k$-ordered, then for every
vertex v  we have $indeg(v),
outdeg(v) \geq k-1$.  \end{corr}

\begin{thm}\label{f_k}

For every $k\geq 2$ there exists an infinite family $\mathcal{F}_k$ of 
$(k-1)$-diregular graphs that are $k$-ordered hamiltonian. 
\end{thm}

\begin{proof}

Given $k$, consider the undirected bracelet graphs $G_{k-1, l}$, $l \geq 3$, as defined in
Section 2.   We define the directed graph $\stackrel{\longrightarrow}{G_{k-1, l}}$ on the
same vertex set with parts $V_1, V_2, \ldots, V_l$, each with $k-1$ vertices. Edge 
$\stackrel{\longrightarrow}{uv}$ is in  $\stackrel{\longrightarrow}{G_{k-1, l}}$  if and only if  $u \in V_i$ and $v \in
V_{i+1}$ where indices
are taken modulo $l$. The  graph
$\stackrel{\longrightarrow}{G_{k-1, l}}$ is $(k-1)$-diregular, and we  now prove that it is 
$k$-ordered hamiltonian. 

Consider marked vertices $v_1, v_2, \ldots, v_k$. We will
 construct a hamiltonian cycle containing these 
vertices in this order. 
It is easy to see that there have to be two consecutive vertices among $v_1, v_2, \ldots, v_k$, 
without loss of generality $v_1$ and $v_2$, 
that are in different parts. We can also suppose without loss of generality that $v_1 \in V_1$.  
Write the vertices of  $\stackrel{\longrightarrow}{G_{k-1, l}}$ as a 
grid where the $i^{th}$ column contains
the vertices in the $i^{th}$ part for $i \in [l]$,  and where the first row contains
the two marked vertices $v_1$ and $v_2$, and where the
$j$th row contains $v_{j+1}$ for $j=2, 3, \ldots, k-1$. 
 Then a hamiltonian cycle containing $v_1, v_2, \ldots, v_k$ in this order is as follows. 
Start at $v_1$ going from left to right across the first row, then 
from the rightmost element in the first row go to the leftmost element in the second row, 
then go from left to right across the second row, 
and so on until reaching the lower rightmost vertex, from which we
close off the hamiltonian cycle by going back to $v_1$.  
 
\end{proof}

\section{Connectivity, linkage,  and $k$-edge-orderedness}

Connectivity, linkage and $k$-orderedness appear to be related concepts.  As
noted in \cite{faun},  a $k$-linked graph $G$ is also $k$-ordered. 
Let $f(k)$ be the minimum connectivity of a graph $G$ that implies $k$-orderedness; the existence
of the function $f(k)$ has been shown in \cite{faun}.   By  a
result in \cite{bol}, $22k$-connected graphs are $k$-linked, and thus are
$k$-ordered as well, leading to the upperbound $f(k)\leq 22k$ observed in
\cite{faun}.

It is natural to pose analogous questions for edge-orderedness and directed
graphs.  In this section we consider edge-orderedness, while
in the following sections we consider directed graphs. 

A graph $G$ is said to be \textit{weakly k-linked} if, given $2k$ vertices (not necessarily distinct)
$s_1, s_2, \ldots, s_k, t_1, t_2, \ldots, t_k$, there exist edge-disjoint paths 
from $s_i$ to $t_i$, for $i\in [k]$.

\begin{lemma}
If a graph $G$ is weakly $2k$-linked, then $G$ is $k$-edge-ordered.
\end{lemma}

\begin{proof}

Consider distinct edges $e_1, \ldots, e_k$.  Let $v_i$ and $u_i$
be the end vertices of $e_i$. As $G$ is weakly $2k$-linked, there are edge-disjoint
paths from $v_1$ to $u_1$, $u_1$ to $v_2, \ldots$, $v_k$ to $u_k$, and $u_k$
to $v_1$.  
If for all
$i$ the path we chose between $v_i$ and $u_i$ is the edge $v_iu_i$,
then we are done.  If the path from $v_i$ to $u_i$ is not
the edge $e_i$, but the edge $e_i$ has not been used in any path, then
we can just replace the path from $v_i$ and $u_i$ with the edge $e_i$. On the
other hand, if  the edge $e_i$ has been used in some other path,
then it has been used by exactly one of them, say $p$.  In this case we can
replace the edge $e_i$ in path $p$ by the path that was between $v_i$ and
$u_i$ and we can replace the path from $v_i$ to $u_i$ with the edge $e_i$.
Repeating this process as necessary, we obtain a tour containing  $e_1, \ldots, e_k$
in this order, and thus $D$ is $k$-edge-ordered. 
 \end{proof}

Let $g(k)$ be the minimum edge-connectivity of a graph $G$ that implies $k$-edge-orderedness
of $G$. 

\begin{prop}
The upper bound  
$g(k) \leq  2k+2$ holds.
\end{prop}

\begin{proof}

It is known that $(k+2)$-edge-connectivity  implies weakly $k$-linked
\cite{huck}.  Thus, the statement of the proposition is a corollary of Lemma
6.1 and \cite{huck}.  \end{proof}

It is easy to see that   $g(k)\geq k-1$.

\begin{lemma}
If a graph $G$ is $2k$-edge-ordered, then it is weakly $k$-linked.
\end{lemma}

\begin{proof} 
Consider vertices $s_1, s_2, \ldots, s_k$ and $t_1, t_2, \ldots, t_k$. 
Since $2k$-edge-ordered implies $2k$-edge-connected
\cite{c}, the degree of every vertex is at least $2k$, and therefore  
there exist edges  $s_1s'_1, s_2s'_2, \ldots, s_k s'_k, t_1t'_1, t_2t'_2, \ldots, t_k t'_k$
such that 
\[
\{s'_1,s'_2, \ldots, s'_k, t'_1, t'_2, \ldots, t'_k\} \cap
\{s_1, s_2, \ldots, s_k, t_1, t_2, \ldots, t_k\}=\emptyset.
\] 
 Then, as $G$ is $2k$-edge-ordered, it readily follows that it is also weakly $k$-linked. 
  \end{proof}

Lemmas 6.1 and 6.3 exhibit a relation between weak linkage and
edge-orderedness, but tightness of these lemmas
remains uncertain.

\section{Connectivity, diameter and orderedness}

In the previous section we considered the relation between connectivity and
edge-orderedness in undirected graphs. In this section we pose the
question: What connectivity implies $k$-(edge-)orderedness in directed graphs?
 
 A digraph $D$ is  \textit{$k$-connected} if 
every vertex cut has at least $k$ vertices. The minimum size of a vertex cut is the 
  \textit{connectivity} of $D$.  
For $S, T \subset V(D)$, let $[S, T]$ be the set of edges from $S$ to $T$. An 
\textit{edge cut} is  the set $[S, \bar{S}]$ for some nonempty $S\subset V(D)$. 
A digraph $D$ is  \textit{$k$-edge-connected} if 
every edge cut has at least $k$ edges. The minimum size of an edge cut is the 
 \textit{edge-connectivity} of $D$ (\cite{west}, Section 4).  
  
It is not immediately clear that any connectivity in directed graphs
implies any orderedness property. Indeed, it has been shown by
Thomassen \cite{th} that for every natural number $k$ there exists a strongly
$k$-connected digraph $D_k$ containing two vertices not lying on a cycle. This implies that
there is no connectivity that would guarantee even the  $2$-(edge-)orderedness of a
directed graph.  When the diameter is small, however, we can prove some positive results.

\begin{thm}
If a digraph $D$ is $g(k)$-edge-connected with diameter $d$ and 
$g(k)\geq (2k-1)\lceil\frac{d}{2}\rceil+1$, then $D$ 
is $k$-edge-ordered.  
\end{thm}

\begin{proof}

Suppose that a digraph $D$ is $g(k)$-edge-connected and it has diameter $d$. Let
$e_1=\stackrel{\longrightarrow}{v_1 u_1}, \ldots, e_k=\stackrel{\longrightarrow}{v_k u_k}$.  As $g(k)\geq 1$, there
exists a directed path between $v_1$ and $u_1$ not longer than $d$.  Delete
the edges of this path from $D$, which decreases  the edge-connectivity
by at most
$\lceil\frac{d}{2}\rceil$. 
Note that $D$ still has edge-connectivity at least $1$, and indeed the edge-connectivity will allow 
us to  repeat the same argument $2k$ times for
paths from $v_1$ to $u_1$, $u_1$ to $v_2, \ldots$, $v_k$ to $u_k$, and $u_k$
to $v_1$.  At the last step the connectivity will be greater than or equal to
$g(k)-(2k-1)\lceil\frac{d}{2}\rceil\geq 1$. Therefore, if $g(k)\geq
(2k-1)\lceil\frac{d}{2}\rceil+1$, then we can obtain an oriented tour
through the vertices $v_1, u_1, \ldots, v_k, u_k$ in this order. 

If for all
$i$ the directed path we chose between $v_i$ and $u_i$ is the edge $\stackrel{\longrightarrow}{v_iu_i}$,
then we are done.  If the path from $v_i$ to $u_i$ is not
the edge $\stackrel{\longrightarrow}{v_iu_i}$, but the edge $\stackrel{\longrightarrow}{v_iu_i}$ has not been used in any path, then
we can just replace the path from $v_i$ and $u_i$ with the edge $\stackrel{\longrightarrow}{v_iu_i}$. On the
other hand, if  the edge $\stackrel{\longrightarrow}{v_iu_i}$ has been used in some other path,
then it has been used by exactly one of them, say $p$.  In this case we can
replace the edge $\stackrel{\longrightarrow}{v_iu_i}$ in path $p$ by the path that was between $v_i$ and
$u_i$ and we can replace the path from $v_i$ to $u_i$ with the edge $\stackrel{\longrightarrow}{v_iu_i}$.
Repeating this process as necessary, we obtain a tour containing  $e_1, \ldots, e_k$
in this order, and thus $D$ is $k$-edge-ordered. 

\end{proof}

\begin{thm}
If a digraph $D$ is $g(k)$-connected with diameter $d\geq 1$, where 
$g(k)\geq (k-1)d$, then $D$ 
is $k$-ordered.  
\end{thm}

\begin{proof} Suppose that a digraph $D$ is $g(k)$-connected and that it has
diameter $d$. Choose marked vertices $v_1, v_2, \ldots,v_k$.
  
Since $g(k)\geq k-1$, removing vertices $v_3, v_4, \ldots, v_k$ would not disconnect $D$, 
thus there exists a directed path from $v_1$ to $v_2$ not containing  $v_3, v_4, \ldots, v_k$
 and furthermore this path has length at most $d$. 
 Delete the vertices of the path from $v_1$ to $v_2$, except $v_1$ and $v_2$. The connectivity
decreases by at most $d-1$, and we can repeat the same process $k$ times for
paths from $v_1$ to $v_2$, $v_2$ to $v_3, \ldots$, $v_k$ to $v_1$.  At the
last step the connectivity will be greater than or equal to
$g(k)-(k-1)(d-1)\geq k-1$, as required. Therefore, if $g(k)\geq (k-1)d$, then we can obtain a
cycle through the vertices $v_1, v_2,  \ldots, v_k$ in this order. This shows
$k$-orderedness.

\end{proof}

Using analogous methods, similar results can be obtained  for undirected graphs. 
Indeed, we get a bound $g(k)\geq (2k-1)d+1$ for a statement analogous to Theorem 7.1 
for undirected graphs, and the same bound in the analogue of Theorem 7.2.

\section{Conclusion}

We conclude by giving an overview of questions motivated by
this paper. In Section 2 we constructed an infinite family of
$2k$-regular $(2k+1)$-ordered hamiltonian bracelet graphs, and in Section 3
we showed that there are no $2k$-ordered bracelet graphs with minimum degree
$2k-1$ and maximum degree less than $2k+2$. We constructed an infinite
family of $2k$-ordered  braceket graphs with minimum degree $2k-1$ and maximum
degree $2k+2$. The following question, however, remains open.  

\begin{q} Is there an infinite family of $2k$-ordered  hamiltonian
bracelet graphs with minimum degree $2k-1$ and maximum degree $2k+2$ 
for all $k\geq 2$?
\end{q}

In Theorem 2.3 we gave a sufficient condition for a bracelet graph to be
$(2k+1)$-ordered.  Note that  Theorem 2.3 only applies when each part
has at least $k$ vertices.  It is not hard to see from connectivity properties
that any $(2k+1)$-ordered bracelet graph with at least $5$ parts has at most two
parts with fewer than $k$ vertices, and if it has two such parts then they
must be adjacent. 

\begin{q} What are the necessary   and sufficient conditions for a bracelet graph
to be $(2k+1)$-ordered?  \end{q}
 
 Naturally, one can also pose this question  for $k$-ordered graphs in general. 
  
 \vskip11pt
In Section 3 we showed that a $2k$-ordered graph has diameter at most 
$\lfloor \frac{n-3}{2k} \rfloor+2$, where $n$ is the number of the vertices of the graph.
We have also shown that there exists an infinite family of $(2k+1)$-ordered 
hamiltonian graphs that have diameter $\lfloor \frac{n-3}{2k} \rfloor+1$. 
It is natural to pose the following question.

\begin{q}  
Is there a $(2k+1)$-ordered or $2k$-ordered 
graph that has diameter $\lfloor \frac{n-3}{2k} \rfloor+2$?  
\end{q}

Note that one can ask the analogues of these questions for directed graphs.

\section{Acknowledgments} This research was performed at the University of
Minnesota Duluth under the supervision of Professor Joseph A. Gallian. The author
would like to thank Professor Gallian for his support and encouragement as
well as Denis Chebikin, Philip Matchett and Melanie Wood for many useful
suggestions. Financial support was provided by the Massachusetts Institute of Technology and 
 the University of Minnesota Duluth.

%\fx{It is a good idea to spellcheck your paper before sending it off.  The
%program I would recommend is MicroSpell, which is available for free for a 45
%day trial period.  See {\tt http://www.microspell.com/} to download the
%software.}

\end{document}